\documentclass[11pt,eqno]{article}
\usepackage{amssymb}
\usepackage{amsmath}
\usepackage{array}
\usepackage{graphicx}
\usepackage{cases}
\usepackage{color}
\usepackage[colorlinks,
            linkcolor=red,
            anchorcolor=green,
            citecolor=blue
            ]{hyperref}
\usepackage{multirow}
\usepackage{picinpar}
\usepackage{fourier}
\makeatletter \renewcommand\@biblabel[1]{[#1]} \makeatother

\begin{document}

\newcommand{\defi}{\stackrel{\Delta}{=}}
\newcommand{\qed}{\hphantom{.}\hfill $\Box$\medbreak}
\newcommand{\A}{{\cal A}}
\newcommand{\B}{{\cal B}}
\newcommand{\U}{{\cal U}}
\newcommand{\G}{{\cal G}}
\newcommand{\cZ}{{\cal Z}}
\newcommand{\proof}{\noindent{\bf Proof \ }}
\newcommand\one{\hbox{1\kern-2.4pt l }}
\newcommand{\Item}{\refstepcounter{Ictr}\item[(\theIctr)]}
\newcommand{\QQ}{\hphantom{MMMMMMM}}

\newtheorem{Theorem}{Theorem}[section]
\newtheorem{Lemma}[Theorem]{Lemma}
\newtheorem{Corollary}[Theorem]{Corollary}
\newtheorem{Remark}[Theorem]{Remark}
\newtheorem{Proposition}[Theorem]{Proposition}
\newtheorem{Assumption}[Theorem]{Assumption}
\newtheorem{Definition}[Theorem]{Definition}
\newtheorem{Construction}[Theorem]{Construction}
\newtheorem{Exa}[Theorem]{Example}
\newcounter{claim_nb}[Theorem]
\setcounter{claim_nb}{0}
\newtheorem{claim}[claim_nb]{Claim}
\newenvironment{cproof}
{\begin{proof}
 [Proof.]
 \vspace{-3.2\parsep}}
{\renewcommand{\qed}{\hfill $\Diamond$} \end{proof}}
\newcommand{\erhao}{\fontsize{21pt}{\baselineskip}\selectfont}
\newcommand{\xiaoerhao}{\fontsize{18pt}{\baselineskip}\selectfont}
\newcommand{\sanhao}{\fontsize{15.75pt}{\baselineskip}\selectfont}
\newcommand{\sihao}{\fontsize{14pt}{\baselineskip}\selectfont}
\newcommand{\xiaosihao}{\fontsize{12pt}{\baselineskip}\selectfont}
\newcommand{\wuhao}{\fontsize{10.5pt}{\baselineskip}\selectfont}
\newcommand{\xiaowuhao}{\fontsize{9pt}{\baselineskip}\selectfont}
\newcommand{\liuhao}{\fontsize{7.875pt}{\baselineskip}\selectfont}
\newcommand{\qihao}{\fontsize{5.25pt}{\baselineskip}\selectfont}

\newcounter{Ictr}

\renewcommand{\theequation}{
\arabic{equation}}
\renewcommand{\thefootnote}{\fnsymbol{footnote}}

\def\A{\mathcal{A}}

\def\C{\mathcal{C}}

\def\V{\mathcal{V}}

\def\I{\mathcal{I}}

\def\Y{\mathcal{Y}}

\def\X{\mathcal{X}}

\def\J{\mathcal{J}}

\def\Q{\mathcal{Q}}

\def\W{\mathcal{W}}

\def\S{\mathcal{S}}

\def\T{\mathcal{T}}

\def\L{\mathcal{L}}

\def\M{\mathcal{M}}

\def\N{\mathcal{N}}
\def\R{\mathbb{R}}
\def\H{\mathbb{H}}


\begin{center}
\topskip2cm
\LARGE{\bf Sparse and Low-Rank  Covariance Matrices Estimation }
\end{center}
\begin{center}
\renewcommand{\thefootnote}{\fnsymbol{footnote}}Shenglong Zhou$^{\textcolor[rgb]{0.98,0.00,0.00}{*}}$,~Naihua Xiu$^{\textcolor[rgb]{0.98,0.00,0.00}{*}}$, ~Ziyan Luo$^{\textcolor[rgb]{0.98,0.00,0.00}{+}}$,~Lingchen Kong\footnote{Aug 7, 2014. E-mail: longnan\_zsl@163.com, nhxiu@bjtu.edu.cn, zyluo@bjtu.edu.cn ,konglchen@126.com}\\
$\textcolor[rgb]{0.98,0.00,0.00}{^*}~$Department of Applied Mathematics\\
${\textcolor[rgb]{0.98,0.00,0.00}{^+}}$~State Key Laboratory of Rail Traffic Control and Safety\\
${\textcolor[rgb]{0.98,0.00,0.00}{^{\underset{+}{*}}}}$~Beijing Jiaotong University, Beijing 100044, P. R. China\\

\end{center}

\vskip12pt

\begin{abstract}
This paper aims at achieving a simultaneously sparse and low-rank estimator from the semidefinite population covariance matrices.
We first benefit from a convex optimization which develops $\ell_1$-norm penalty to encourage the sparsity and nuclear norm to favor the low-rank property. For the proposed estimator, we then prove that with large probability, the Frobenious norm of the estimation rate can be of order $\mathcal{O}\left(\sqrt{(s\log r)/n}\right)$ under a mild case, where $s$ and $r$ denote the number of sparse entries and the rank of the population covariance respectively, $n$ notes the sample capacity. Finally an efficient alternating direction method of multipliers with global convergence is proposed to tackle this problem, and meantime merits of the approach are also illustrated by practicing numerical simulations.
\end{abstract}
\section{Introduction}\label{sec1}
Estimation of population covariance matrices from samples of multivariate data has draw many attentions in the last decade owing to its fundamental importance in multivariate analysis. With dramatic advances in technology in recent years, various research fields, such as genetic data, brain imaging, spectroscopic imaging, climate data and so on, have been used to deal with massive high-dimensional data sets, whose sample sizes can be very small relative to dimension. In such settings, the standard and the most usual sample covariance matrices often performs poorly \cite{BLa,BLb,J}.

Fortunately, regularization as a class of new methods to estimate covariance matrices has recently emerged to overcome those shortages of using traditional sample covariance matrices. These methods encompass several specified forms, banding \cite{BLa,CZa,WP}, tapering \cite{CZZ,FB} and thresholding
\cite{BLb,CY,E,RLZ} for instance.

Moreover, there are many cases where the model is known to be structured in several ways at the same time. In recent years, one of research contents is to estimate a covariance matrix possessing both sparsity and positive definiteness. For instance, Rothman \cite{R} gave the following model:
\begin{eqnarray}\label{Rest}\hat{\Sigma}_+=\underset{\Sigma\succ0}{\textrm{argmin}}~\frac{1}{2}\left\| \Sigma-\Sigma_{n}\right\|_{F}^{2}-\tau\log\det(\Sigma)+\lambda\left\|\Sigma^-\right\|_{1},\end{eqnarray}
where $\Sigma_n$ is the sample covariance matrix, $\|\cdot\|_F$ is the Frobenious norm, $\|\cdot\|_1$ is the element-wise $\ell_1$-norm and $\Sigma^-:=\Sigma-\textmd{Diag}(\Sigma)$. From the optimization viewpoint, (\ref{Rest}) is similar to the graphical lasso criterion \cite{FHT} which also has a log-determinant part and the element-wise $\ell_1$-penalty. Rothman \cite{R} derived an iterative procedure to solve (\ref{Rest}). While Xue, Ma and Zou \cite{XMZ} omitted the log-determinant part and considered the positive definite constraint $\{\Sigma\succeq\varepsilon \textbf{\emph{I}}\}$ for some arbitrarily small $\varepsilon>0$:
\begin{eqnarray}\label{Xest}\hat{\Sigma}_{\varepsilon}=\underset{\Sigma\succeq\varepsilon \textbf{\emph{I}}}{\textrm{argmin}}~\frac{1}{2}\left\| \Sigma-\Sigma_{n}\right\|_{F}^{2}+\lambda\left\|\Sigma^-\right\|_{1},\end{eqnarray}
They utilized an efficient alternating direction method (ADM) to solve the challenging problem (\ref{Xest}) and established its convergence properties.

Most of the literatures, e.g.,\cite{LWZ,R,XMZ}, required the population covariance matrices being positive definite, and thus there is no essence of pursuing the low-rank of the estimator. By contrast, newly appeared research topic is to consider simultaneously the sparsity and low-rank of a structured model, which implies that the population covariance matrices are no longer restricted to the positive definite matrix cone and can be relaxed to the positive semidefinite cone. In addition, the models with structure of being simultaneously the sparsity and low-rank are widely applied into practice, such as sparse signal recovery from quadratic measurements and sparse phase retrieval, see \cite{O} for example. Moreover, Richard et al. \cite{RSV} showed that both sparse and low-rank model can be derived in covariance matrix when the random variables are highly correlated in groups, which means this covariance matrix has a block diagonal structure.

With stimulations of those ideas, we construct the following convex model  encompassing the $\ell_1$-norm and nuclear norm  for estimating the covariance matrix:
\begin{eqnarray}\label{est}\hat{\Sigma}=\underset{\Sigma\succeq0}{\textrm{argmin}}~\frac{1}{2}\left\| \Sigma-\Sigma_{n}\right\|_{F}^{2}+\lambda\left\|\Sigma\right\|_{1}+\tau\left\|\Sigma\right\|_{*},\end{eqnarray}
where $\lambda\geq0, \tau\geq0$ are tuning parameters. The $\ell_1$-norm penalty $\left\|\Sigma\right\|_{1}=\sum_{i,j}|\sigma_{ij}|$ is also called lasso-type penalty and is used to encourage sparse solutions. The nuclear norm , $\left\|\Sigma\right\|_{*}=\sum_{i}|\lambda_{i}(\Sigma)|$ with $\lambda_{i}(\Sigma)$ being the eigenvalue of $\Sigma$, is the trace norm when $\Sigma\succeq0$ and ensures low-rank solutions of (\ref{est}). Here we inroduce the \emph{approximate rank} to interpret the low-rank, which is defined as $ar(A):=ar\in\{1,2,\cdots,p\}$ being the smallest number such that
\begin{eqnarray}\label{apr}\frac{\sigma_{ar+1}(A)}{\sigma_1(A)}\leq \gamma,\end{eqnarray}
where $\sigma_i(A)(i=1,2,\cdots,p)$ are the singular value of $A$ with $\sigma_1(A)\geq\sigma_2(A)\geq\cdots\geq\sigma_p(A)$, and $\gamma>0$ could be chosen based on the needs, throughout our paper we fix $\gamma=0.001$ for simplicity.

The contributions of this paper mainly center on two aspects. For one thing, being different from \cite{O,RSV}, we establish the theoretical statistical theory under different assumptions rather than giving the generalized error bound of the estimation. Especially, we acquire the estimation rate $\mathcal{O}(\sqrt{(s\log r)/n})$ under the Frobenious norm error, which improves the optimal rate $\mathcal{O}(\sqrt{(s\log p)/n})$ where the low-rank property of the estimator does not be considered \cite{CZb,R,XMZ} and $p$ is the samples' dimension with $p>\max\{n,r\}$. For another, we take advantage of the alternating direction method of multipliers (ADMM), also can be seen in \cite{XMZ,Y}, to combat our problem (\ref{est}).

The organization of this paper is as follows. In Section \ref{sec2} we will present some theoretical properties of the estimator derived by the proposed model (\ref{est}). After that the alternating direction method of multipliers (ADMM) is going to be introduced to combat the problem, and numerical experiments are projected to show the performance of this method in Sections \ref{sec3} and \ref{sec4} respectively. We make a conclusion in the last section.
\section{A Sparse and Low-Rank Covariance Estimator }\label{sec2}
Before the main part, we hereafter introduce some notations. $\mathbb{E}(X)$ and $\mathbb{P}(A)$ denote the expectation of $X$ and the probability of the incident $A$ occurring respectively. $\textmd{Card}(S)$ is the number of entries of the set $S$. Normal distribution with mean $\mu$ and covariance $\Sigma$ is written as $N\left(\mu,\Sigma\right)$. Say $Y_n=\mathcal{O}_P(1)$ if for every $\varepsilon> 0$, there is a $C > 0$ such that
$\mathbb{P}\{|Y_n|>C\}<\varepsilon$ for all $n\geq n_0(\varepsilon)$, and say $Y_n=\mathcal{O}_P(a_n)$ if $Y_n/a_n=\mathcal{O}_P(1)$. If there are two constants $C_1\leq C_2$ such that $C_1\leq X_n/Y_n\leq C_2$, we write as $X_n\asymp Y_n$.

For given observed independently and identically distributed (i.i.d. for short) $p$-variate random variables $X_{1},\cdots,$ $X_{n}$ with covariance matrix $\Sigma_{0}$ and $p>n$, the goal is to estimate the unknown matrix $\Sigma_{0}$ based on the sample $\{X_{l}: l=1, \cdots, n\}$. This problem is called covariance matrix estimation which is of fundamental importance in multivariate analysis.

Given a random sample $\left\{X_{1},\cdots,X_{n}\right\}$ from $\mathbb{E}(X)=0$ (without loss of generality) and a population covariance matrix $\Sigma_0=\left(\sigma_{0ij}\right)_{1\leq i,j\leq p}=\mathbb{E}\left(XX^{\top}\right)$, the sample covariance matrix is
\[\Sigma_n=\left(\sigma_{nij}\right)_{1\leq i,j\leq p}=\frac{1}{n-1}\sum_{l=1}^{n}\left(X_{l}-\bar{X}\right)\left(X_{l}-\bar{X}\right)^{\top},\]
where $\bar{X}=\frac{1}{n}\sum_{l=1}^{n}X_{l}$. Denote $S$  the support set of the population covariance matrix $\Sigma_0$ as \[S=\left\{(i,j):\sigma_{0ij}\neq 0\right\},~s=\textmd{Card} (S),~r=\textmd{rank}(\Sigma_0).\]

\begin{Assumption}\label{Ass1} For all $p$, $0\leq\lambda_{\emph{min}}\left(\Sigma_{0}\right)\leq\lambda_{\emph{max}}\left(\Sigma_{0}\right)\leq\bar{\lambda} <\infty$, where $\bar{\lambda} $ is a constant.\end{Assumption}
Assumption \ref{Ass1} is a common used condition in covariance matrix estimation, on which a useful lemma based is recalled here for the sequel analysis. One can also refer it in \cite{BLa}.

\begin{Lemma}\label{ZZ} Let~~$X_l\in\mathbb{R}^p, l=1,\cdots,n$ be i.i.d. $N\left(0, \Sigma_0\right)$ and Assumption \ref{Ass1} holds, \emph{(}i.e., for all $p$, $\lambda_{\emph{max}}\left(\Sigma_{0}\right)\leq\bar{\lambda} <\infty$\emph{)}.  Then, if~~$\Sigma_0=\left(\sigma_{0ij}\right)_{1\leq i,j\leq p}$,
\begin{eqnarray}\label{ZZZ}\mathbb{P}\left\{\left|\sum_{l=1}^{n}\left(X_{li}X_{lj}-\sigma_{0ij}\right)\right|\geq n\nu\right\}\leq C_{1}\emph{\textmd{exp}}\left\{-C_2n\nu^2\right\},~for~|\nu|\leq\delta\end{eqnarray}
where constants $C_1, C_2$ and $\delta$ depend on $\bar{\lambda}$ only.\end{Lemma}

Another lemma which plays an important role in our main results is stated below.
\begin{Lemma}\label{IM} Suppose that Assumption \ref{Ass1} holds, $\lambda\leq\frac{\varepsilon}{16\sqrt{s}}, \tau\leq\frac{\varepsilon}{8\sqrt{r}}$. Then for $\varepsilon>0$ sufficiently small,
\begin{eqnarray}\label{max}\underset{1\leq i,j\leq p}{\textrm{\emph{max}}}\left|\sigma_{0ij}-\sigma_{nij}\right|\leq \lambda~~implies~~ \|\hat{\Sigma}-\Sigma_0\|_F\leq\varepsilon,\end{eqnarray}
where $\hat{\Sigma}$ is defined as $(\ref{est})$.\end{Lemma}
\noindent\textbf{Proof}~~First make the eigenvalue decomposition of $\Sigma_0$ ($r=\textmd{rank}(\Sigma_0)$) as
$$\Sigma_0=U
\Lambda_{\Sigma_0}U^\top\equiv U
\left(
              \begin{array}{cc}
                \textmd{Diag}(\lambda(\Sigma_0)) & ~ \\
                ~ & 0 \\
              \end{array}
            \right)U^\top,$$
where $U\in\mathbb{R}^{p\times p}$ with $UU^\top=U^\top U=\textbf{I}_{p}$ is the matrix composed of eigenvectors, $\textmd{Diag}(\lambda(\Sigma_0))\in\mathbb{R}^{r\times r}$ is a diagonal matrix generated by eigenvalues with $\lambda(\Sigma_0)=(\lambda_1(\Sigma_0),\cdots,\lambda_r(\Sigma_0))^\top$ and $\lambda_j(\Sigma_0)>0, j=1,\cdots,r$.\\
\noindent By denoting $\Sigma :=U\Delta U^{\top}+\Sigma_0$ with $\Delta=\left(
                                                                       \begin{array}{cc}
                                                                         \Delta_1 & 0 \\
                                                                         0 & 0 \\
                                                                       \end{array}
                                                                     \right)
$ and $\Delta_1\in\mathbb{R}^{r\times r}$,  which implies $\Delta=U^{\top}\Sigma U-\Lambda_{\Sigma_0}$, we  consider the model
\begin{eqnarray}\label{delta}\hat{\Delta}&:=&\underset{\Delta+\Lambda_{\Sigma_0}\succeq0}{\textmd{argmin}}~~ F(\Delta)=\underset{\Delta_1+\textmd{Diag}(\lambda(\Sigma_0))\succeq0}{\textmd{argmin}}~~ F(\Delta)\\
&\equiv&\underset{\Delta+\Lambda_{\Sigma_0}\succeq0}{\textmd{argmin}}~~\frac{1}{2}\left\| U\Delta U^{\top}+\Sigma_0-\Sigma_n\right\|_{F}^{2}+\lambda\left\|U\Delta U^{\top}+\Sigma_0\right\|_{1}+\tau\left\|U\Delta U^{\top}+\Sigma_0\right\|_{*}.\nonumber\end{eqnarray}
Clearly, from (\ref{est}) we have $\hat{\Sigma} =U\hat{\Delta} U^{\top}+\Sigma_0$ which implies $\hat{\Delta}=U^{\top}\hat{\Sigma} U-\Lambda_{\Sigma_0}$. For a given $\varepsilon>0$  sufficiently small and any $\|\Delta\|_F=\varepsilon$ (i.e., $\|\Delta_1\|_F=\varepsilon$), we compute
\begin{eqnarray}
F(\Delta)-F(0)&=&\frac{1}{2}\left\| U\Delta U^{\top}+\Sigma_0-\Sigma_n\right\|_{F}^{2}+\lambda\left\|U\Delta U^{\top}+\Sigma_0\right\|_{1}+\tau\left\|U\Delta U^{\top}+\Sigma_0\right\|_{*}\nonumber\\
&&-\frac{1}{2}\left\|\Sigma_0-\Sigma_n\right\|_{F}^{2}+\lambda\left\|\Sigma_0\right\|_{1}+\tau\left\|\Sigma_0\right\|_{*}\nonumber\\
&=&\frac{1}{2}\| \Delta \|_{F}^{2}+\left\langle U\Delta U^{\top}, \Sigma_0-\Sigma_n\right\rangle+\lambda\left(\|U\Delta U^{\top}+\Sigma_0\|_{1}-\|\Sigma_0\|_{1}\right)\nonumber\\
&&+\tau\left(\|U\Delta U^{\top}+\Sigma_0\|_{*}-\|\Sigma_0\|_{*}\right)\nonumber\\
&\equiv&\frac{1}{2}\| \Delta \|_{F}^{2}+\textmd{I}+\textmd{II}+\textmd{III}.\nonumber\end{eqnarray}
For convenience we denote $\overline{\Delta}=U\Delta U^{\top}$. Then for I,  it holds
\begin{eqnarray}\textmd{I}=\left\langle \overline{\Delta}, \Sigma_0-\Sigma_n\right\rangle&\geq&-\left|\textmd{tr}\left(\overline{\Delta}\left(\Sigma_n-\Sigma_0\right)\right)\right|\nonumber\\
&=&-\Big|\sum_{i,j}\overline{\Delta}_{ij}\left(\sigma_{0ij}-\sigma_{nij}\right)\Big|
\geq-\underset{1\leq i,j\leq p}{\textrm{max}}\left|\sigma_{0ij}-\sigma_{nij}\right|\|\overline{\Delta}\|_1,\nonumber\end{eqnarray}
For II, we obtain by noting $S=\left\{(i,j): \sigma_{0ij}\neq 0\right\}$ and $s=\textmd{Card}(S)$ that,
\begin{eqnarray}\textmd{II}
&=&\lambda\left(\|\Sigma_0+\overline{\Delta}\|_{1}-\left\|\Sigma_0\right\|_{1}\right)\nonumber\\
&=&\lambda\left(\|\Sigma_{0S}+\overline{\Delta}_S\|_{1}+\|\overline{\Delta}_{S^C}\|_{1}-\left\|\Sigma_{0S}\right\|_{1}\right)\nonumber\\
&\geq&\lambda\Big(\|\Sigma_{0S}+\overline{\Delta}_S\|_{1}+\|\overline{\Delta}_{S^C}\|_{1}-
(\|\Sigma_{0S}+\overline{\Delta}_S\|_{1}+\|\overline{\Delta}_S\|_{1})\Big)\nonumber\\
&=&\lambda\left(\|\overline{\Delta}_{S^C}\|_{1}-\|\overline{\Delta}_S\|_{1}\right).\nonumber\end{eqnarray}
From the H$\ddot{o}$lder Inequality, one can prove that
\begin{eqnarray}\label{HD1}&&\|\overline{\Delta}\|_{*}=\|\Delta\|_{*}=\|\Delta_1\|_{*}=\sum_{i=1}^r\left|\lambda_i(\Delta_1)\right|\leq
\sqrt{r}\|\Delta_1\|_{F}=\sqrt{r}\|\Delta\|_{F},\\
\label{HD2}&&\|\overline{\Delta}_S\|_1\leq\sqrt{s}\|\overline{\Delta}_S\|_F.\end{eqnarray}
For III, combining with (\ref{HD1}) we get that
\begin{eqnarray}\textmd{III}= \tau\left(\|\Sigma_0+\overline{\Delta}\|_{*}-\left\|\Sigma_0\right\|_{*}\right)
\geq-\tau\|\overline{\Delta}\|_{*}=-\tau\|\Delta\|_{*}\geq-\tau\sqrt{r}\|\Delta\|_{F}.\nonumber\end{eqnarray}
Since $\textrm{max}_{i,j}\left|\sigma_{0ij}-\sigma_{nij}\right|\leq \lambda$ and (\ref{HD2}),
\begin{eqnarray}G(\Delta)&\equiv& F(\Delta)-F(0)\nonumber\\
&=&\frac{1}{2}\left\|\Delta\right\|_{F}^{2}+\textmd{I}+\textmd{II}+\textmd{III}\nonumber\\
&\geq&\frac{1}{2}\left\|\Delta\right\|_{F}^{2}-\lambda\|\overline{\Delta}\|_{1}+\lambda\left(\|\overline{\Delta}_{S^C}\|_{1}-\|\overline{\Delta}_S\|_{1}\right)
-\tau\sqrt{r}\|\Delta\|_{F}\nonumber\\
&=&\frac{1}{2}\left\|\Delta\right\|_{F}^{2}-\lambda\left(\|\overline{\Delta}_S\|_{1}+\|\overline{\Delta}_{S^C}\|_{1}
-(\|\overline{\Delta}_{S^C}\|_{1}-\|\overline{\Delta}_S\|_{1})\right)
-\tau\sqrt{r}\|\Delta\|_{F}\nonumber\\
&=&\frac{1}{2}\left\|\Delta\right\|_{F}^{2}-2\lambda\|\overline{\Delta}_S\|_{1}-\tau\sqrt{r}\|\Delta\|_{F}\nonumber\\
&\geq&\frac{1}{2}\left\|\Delta\right\|_{F}^{2}-2\lambda\sqrt{s}\|\overline{\Delta}_S\|_{F}-\tau\sqrt{r}\|\Delta\|_{F}\nonumber\\
&\geq&\frac{1}{2}\left\|\Delta\right\|_{F}^{2}-2\lambda\sqrt{s}\|\overline{\Delta}\|_{F}-\tau\sqrt{r}\|\Delta\|_{F}\nonumber\\
&=&\frac{1}{2}\left\|\Delta\right\|_{F}^{2}-2\lambda\sqrt{s}\|\Delta\|_{F}-\tau\sqrt{r}\|\Delta\|_{F}.\nonumber\end{eqnarray}
Therefore, by $\lambda\leq\frac{\varepsilon}{16\sqrt{s}}, \tau\leq\frac{\varepsilon}{8\sqrt{r}}$
$$G(\Delta)\geq\frac{1}{2}\left\|\Delta\right\|_{F}^{2}-2\lambda\sqrt{s}\|\Delta_S\|_{F}-\tau\sqrt{r}\left\|\Delta\right\|_{F}
\geq\frac{\varepsilon^{2}}{2}-\frac{\varepsilon^{2}}{8}-\frac{\varepsilon^{2}}{8}=\frac{\varepsilon^{2}}{4}>0.$$
Hence we prove that if $\lambda\leq\frac{\varepsilon}{16\sqrt{s}}, \tau\leq\frac{\varepsilon}{8\sqrt{r}}$ and $\textrm{max}_{i,j}\left|\sigma_{0ij}-\sigma_{nij}\right|\leq \lambda$, for any $\|\Delta\|_F=\varepsilon$, it holds
$$G(\Delta)\equiv F(\Delta)-F(0)>0.$$
In addition, from (\ref{delta}), we have
\begin{eqnarray}\label{delta1}\hat{\Delta}=\underset{\Delta+\textmd{Diag}(\lambda_{\Sigma_0})\succeq0}{\textmd{argmin}}~~F(\Delta)-F(0) ~=\underset{\Delta+\textmd{Diag}(\lambda_{\Sigma_0})\succeq0}{\textmd{argmin}}~~G(\Delta),\end{eqnarray}
which implies that $\|\hat{\Delta}\|_F\leq\varepsilon$. Otherwise, we suppose $\|\hat{\Delta}\|_F>\varepsilon$, then $\|\varepsilon\hat{\Delta}/\|\hat{\Delta}\|_F\|_F=\varepsilon$. Since for any $\|\Delta\|_F=\varepsilon$, it follows $G(\Delta)>0=G(0)$ which is contradicted with the fact $G(\cdot)$ is a convex function and $G(\hat{\Delta})\leq G(0)=0$, because
\begin{eqnarray}0<G\left(\frac{\varepsilon}{\|\hat{\Delta}\|_F}\hat{\Delta}\right)&=&
G\left(\frac{\varepsilon}{\|\hat{\Delta}\|_F}\hat{\Delta}+\left(1-\frac{\varepsilon}{\|\hat{\Delta}\|_F}\right)0\right)\nonumber\\
&\leq&
\frac{\varepsilon}{\|\hat{\Delta}\|_F}G(\hat{\Delta})+\left(1-\frac{\varepsilon}{\|\hat{\Delta}\|_F}\right)G(0)\leq0.\nonumber\end{eqnarray}
Finally $\|\hat{\Delta}\|_F\leq\varepsilon$ indicates that  $\|\hat{\Sigma}-\Sigma_0\|_F\leq\varepsilon$ due to $\|\hat{\Delta}\|_F=\|U\hat{\Delta}U^{\top}\|_F=\|\hat{\Sigma}-\Sigma_0\|_F$. Hence the desired result is obtained.{\qed}

Then in order to acquiring the rate of the estimation, the two following commonly used assumptions are needed to introduced, and also can be seen \cite{BLa,R}. Assumption \ref{Ass2} holds, for example, if $X_{l1}, \cdots, X_{lp}$ are Gaussian.
\begin{Assumption}\label{Ass2} $\mathbb{E}\left\{\emph{\textrm{exp}}(tX_{lj}^2)\right\}\leq C_{1} <\infty$ hold for all $j =1,\cdots, p$ and $0<|t|<t_0$, where $t_0$ and $C_{1}$  are two constants.\end{Assumption}
\begin{Assumption}\label{Ass3} $\mathbb{E}\left\{\left|X_{lj}\right|^{2\alpha}\right\}\leq C_2 <\infty$ hold for all $j=1,\cdots, p$, where some $\alpha\geq2$ and $C_{2}$ is a constant.\end{Assumption}
Built on the two assumptions, we give our main results with regard to rates of the estimator of (\ref{est}).
\begin{Theorem}\label{Th1} For some $\delta\in[1,2)$, let $K_1$ be a sufficiently large constant and suppose Assumptions \ref{Ass1} and \ref{Ass2} hold. If $\lambda=K_1\sqrt{\frac{\log r}{n}+\frac{\log p}{n K_1^\delta}}, \tau=\mathcal{O}_P\left(\sqrt{\frac{s\log r}{nr}+\frac{s\log p}{nrK_1^\delta}}\right)$ and $s\log r+(s\log p)/K_1^\delta=o(n)$, then
\begin{eqnarray}\label{rate1}\|\hat{\Sigma}-\Sigma_0\|_F=\mathcal{O}_P\left(\sqrt{\frac{s\log r}{n}+\frac{s\log p}{nK_1^\delta}}\right).\end{eqnarray}
\end{Theorem}
\noindent\textbf{Proof}~~Since Assumptions \ref{Ass1} and \ref{Ass2} hold, a fact employed by Rothman et al. \cite{RLZ}
is that for $\nu>0$ sufficiently small,
\begin{eqnarray}\mathbb{P}\left\{\underset{1\leq i,j\leq n}{\textrm{max}}\left|\sigma_{0ij}-\sigma_{nij}\right|\geq \nu\right\}\leq C_{3}p^2\textrm{exp}\left\{-C_4n\nu^2\right\},\nonumber\end{eqnarray}
where $C_{3}$ and $C_{4}$ are some constants. We then apply the bound $s\log r+(s\log p)/K_1^\delta=o(n)$ and Lemma \ref{IM} with $$\varepsilon=16K_1\sqrt{\frac{s\log r}{n}+\frac{s\log p}{nK_1^\delta}}, ~~ \lambda=K_1\sqrt{\frac{\log r}{n}+\frac{\log p}{n K_1^\delta}}$$ to obtain that
\begin{eqnarray}\mathbb{P}\left\{\|\hat{\Sigma}-\Sigma_0\|_F\leq\varepsilon\right\}
\geq\mathbb{P}\left\{\underset{1\leq i,j\leq n}{\textrm{max}}\left|\sigma_{0ij}-\sigma_{nij}\right|\leq\lambda\right\}\geq 1-C_3p ^{2-C_4K_1^{2-\delta}}r^{-C_4K_1^2}.\nonumber\end{eqnarray}
Evidently, $1-C_3p ^{2-C_4K_1^{2-\delta}}r^{-C_4K_1^2}$ can be arbitrarily close to one by choosing $K_1$ sufficiently large.\qed
\begin{Corollary}\label{Co1} For some $\delta\in[1,2)$, let $K_1$ be a sufficiently large constant such that $K_1^\delta\log r\asymp\log p$, and suppose Assumptions \ref{Ass1}, \ref{Ass2} hold. If $\lambda=K_1\sqrt{\frac{\log r}{n}}, \tau=\mathcal{O}_P\left(\sqrt{\frac{s\log r}{nr}}\right)$ and $s\log r=o(n)$, then
\begin{eqnarray}\label{rate11}\|\hat{\Sigma}-\Sigma_0\|_F=\mathcal{O}_P\left(\sqrt{\frac{s\log r}{n}}\right).\end{eqnarray}
\end{Corollary}
\begin{Remark}Clearly, if the $\lambda_{\min}(\Sigma_0)>0$ in Assumption \ref{Ass1}, the better rate $\mathcal{O}_P\left(\sqrt{\frac{s\log r}{n}}\right)$ would reduce to $\mathcal{O}_P\left(\sqrt{\frac{s\log p}{n}}\right)$. It is worth mentioning that under the the Assumption \ref{Ass1}, the minimax optimal rate of convergence under
the Frobenius norm in Theorem 4 of \cite{CZb} is $\mathcal{O}_P\left(\sqrt{\frac{s\log p}{n}}\right)$ which also has been obtained by \cite{XMZ}. However, to attain the same rate in the presence of the log-determinant barrier term (\ref{Rest}), Rothman \cite{R} instead would require that $\lambda_{\min}$, the minimal eigenvalue of the true covariance matrix, should be bounded away from zero by some positive constant, and also that the
barrier parameter should be bounded by some positive quantity. \cite{XMZ} illustrated this theory requiring a lower bound on $\lambda_{\min}$ is not very appealing.\end{Remark}
\begin{Theorem}\label{Th2} Let $K_2$ be a sufficiently large constant and suppose that Assumptions \ref{Ass1} and \ref{Ass3} hold. If $\lambda=K_2\sqrt{\frac{p^{4/\alpha}}{n}}, \tau=\mathcal{O}_P\left(\sqrt{\frac{sp^{4/\alpha}}{nr}}\right)$ and $sp^{4/\alpha}=o(n)$, then
\begin{eqnarray}\label{rate2}\|\hat{\Sigma}-\Sigma_0\|_F=\mathcal{O}_P\left(\sqrt{\frac{sp^{4/\alpha}}{n}}\right).\end{eqnarray}\end{Theorem}
\noindent\textbf{Proof}~~Since Assumptions \ref{Ass1} and \ref{Ass3} hold, one can modify a result of Bickel \& Levina (2008a)
and show that for $\nu>0$ sufficiently small
\begin{eqnarray}\mathbb{P}\left\{\underset{1\leq i,j\leq n}{\textrm{max}}\left|\sigma_{0ij}-\sigma_{nij}\right|\geq \nu\right\}\leq p^2C_{5}n^{-\alpha/2}\nu^{-\alpha},\nonumber\end{eqnarray}
where $C_5$ is a constant. We then apply the bound $sp^{4/\alpha}=o(n)$ and Lemma \ref{IM} with $\varepsilon=16K_2\sqrt{\frac{sp^{4/\alpha}}{n}}, \lambda=K_2\sqrt{\frac{p^{4/\alpha}}{n}}$ to get
\begin{eqnarray}\mathbb{P}\left\{\left\|\hat{\Sigma}-\Sigma_0\right\|_F\leq\varepsilon\right\}
\geq\mathbb{P}\left\{\underset{1\leq i,j\leq n}{\textrm{max}}\left|\sigma_{0ij}-\sigma_{nij}\right|\leq\lambda\right\}\geq 1-C_6K_2^{-\alpha}.\nonumber\end{eqnarray}
Apparently, the bound $1-C_6K_2^{-\alpha}$ can be arbitrarily close to one by taking $K_2$ sufficiently large.{\qed}

\section{Alternating Direction Method of multipliers}\label{sec3}

In this section, we will construct the alternating direction method of multipliers (ADMM) to solve problem (\ref{est}). By introducing an auxiliary variable $\Gamma$, problem (\ref{est}) can be rewritten as
\begin{eqnarray}\label{est1}
&&\textmd{min}~~\frac{1}{2}\left\|\Sigma-\Sigma_{n}\right\|_{F}^{2}+\lambda\left\|\Sigma\right\|_{1}+\tau\|\Gamma\|_{*},\\
&&\textmd{s.t.}~~~~\Gamma\succeq0,~\Sigma-\Gamma=0.\nonumber\end{eqnarray}
The constraint $\Gamma\succeq0$ can be put into the objective function by using an indicator function:
\begin{numcases}{\mathcal{I}\left(\Gamma\succeq0\right)=}
~~0,~~~~~~~~\Gamma\succeq0\nonumber\\
~+\infty,~~~\textmd{otherwise}. \nonumber\end{numcases}
This leads to the following equivalent reformulation of (\ref{est1}):
\begin{eqnarray}\label{est2}
&&\textmd{min}~~\frac{1}{2}\left\|\Sigma-\Sigma_{n}\right\|_{F}^{2}+\lambda\left\|\Sigma\right\|_{1}+\tau\|\Gamma\|_{*}
+\mathcal{I}\left(\Gamma\succeq0\right),\\
&&\textmd{s.t.}~~~~\Sigma-\Gamma=0.\nonumber\end{eqnarray}
Recently, the alternating direction method of multipliers (ADMM) has been studied extensively for solving (\ref{est2}). A typical iteration of ADMM for solving (\ref{est2}) can be described as
\begin{numcases}{}
\label{AD1}\Gamma^{k+1}:=\textmd{argmin}_{\Gamma}~\mathcal{L_{\mu}}(\Sigma^{k},\Gamma,\Lambda^k)\\
\label{AD2}\Sigma^{k+1}:=\textmd{argmin}_{\Sigma}~\mathcal{L_{\mu}}(\Sigma,\Gamma^{k+1},\Lambda^k)\\
\label{AD3}\Lambda^{k+1}:=\Lambda^k-\frac{1}{\mu}\left(\Gamma^{k+1}-\Sigma^{k+1}\right), \end{numcases}
where the augmented Lagrangian function $\mathcal{L}_\mu(\Sigma,\Gamma,\Lambda)$ is defined as
\begin{eqnarray}
\mathcal{L}_\mu(\Sigma,\Gamma,\Lambda)&:=&\frac{1}{2}\left\|\Sigma-\Sigma_{n}\right\|_{F}^{2}+\lambda\left\|\Sigma\right\|_{1}+\tau\|\Gamma\|_{*}\nonumber\\
\label{ALgr}&&+\mathcal{I}(\Gamma\succeq0)-\left\langle\Lambda,\Gamma- \Sigma\right\rangle+\frac{1}{2\mu}\|\Sigma-\Gamma\|_F^2,\end{eqnarray}
in which $\Lambda$ is the Lagrange multiplier and $\mu>0$ is a penalty parameter.
Note that ADMM (\ref{AD1}-\ref{AD3}) can be written explicitly as\\
\begin{numcases}{}
\label{AD11}\Gamma^{k+1}:=\textmd{argmin}_{\Gamma\succeq0}~\tau\|\Gamma\|_*+
\frac{1}{2\mu}\|\Gamma-(\Sigma^{k}+\mu\Lambda^k)\|_{F}^{2}\\
\label{AD22}\Sigma^{k+1}:=\textmd{argmin}_{\Sigma}~\lambda\left\|\Sigma\right\|_1+
\frac{\mu+1}{2\mu}\|\Sigma-\frac{\mu}{\mu+1}(\Sigma_n+\frac{1}{\mu}\Gamma^{k+1}-\Lambda^k)\|_{F}^{2}\\
\label{AD33}\Lambda^{k+1}:=\Lambda^k-(\Gamma^{k+1}-\Sigma^{k+1})/\mu. \end{numcases}
We now show that the two subproblems (\ref{AD11}) and (\ref{AD22})  can be easily solved. For the subproblem  (\ref{AD11}),
\begin{eqnarray} \Gamma^{k+1}&=&\underset{\Gamma\succeq0}{\textmd{argmin}}~\tau\|\Gamma\|_*+
\frac{1}{2\mu}\|\Gamma-(\Sigma^{k}+\mu\Lambda^k)\|_{F}^{2}\nonumber\\
&=&\underset{\Gamma\succeq0}{\textmd{argmin}}~\tau\langle\Gamma,I\rangle+
\frac{1}{2\mu}\|\Gamma-(\Sigma^{k}+\mu\Lambda^k)\|_{F}^{2}\nonumber\\
&=&\underset{\Gamma\succeq0}{\textmd{argmin}}~\|\Gamma-(\Sigma^{k}+\mu\Lambda^k-\mu\tau I)\|_{F}^{2}\nonumber\\
&=&(\Sigma^{k}+\mu\Lambda^k-\mu\tau I)_+,\end{eqnarray}
where $(X)_+$ denote the projection of a matrix $X$ onto the convex positive semidefinite cone $S^n_+$. Namely $(X)_+=U\textmd{Diag}(\max\{\lambda(X),0\})U^{\top}$, where $X=U\textmd{Diag}(\lambda(X))U^{\top}$ and $\lambda(X)=(\lambda_{1}(X),\lambda_{2}(X),\cdots,\lambda_{p}(X))^{\top}.$

The solution of the second subproblem (\ref{AD22}) is given by the $\ell_1$-shrinkage operation
\begin{eqnarray}\label{Shr1}\Sigma^{k+1}&=&\frac{\mu}{\mu+1}\textmd{Shrink}(\Sigma_n+\frac{1}{\mu}\Gamma^{k+1}-\Lambda^k,\lambda)\nonumber\\
&=&\frac{\mu}{\mu+1}\left(\textmd{max}~\left\{\left|P_{ij}\right|-\lambda,0\right\}\textmd{sign}\left(P_{ij}\right)\right)_{1\leq i,j\leq p},\end{eqnarray}
where $P=\left(P_{ij}\right)_{1\leq i,j\leq p}:=\Sigma_n+\frac{1}{\mu}\Gamma^{k+1}-\Lambda^k$ and $\textmd{sign}(\cdot)$ is a sign function.

Therefore, combining with (\ref{AD33})-(\ref{Shr1}), the whole algorithm is written as follows

\begin{table}[h]\renewcommand{\arraystretch}{1.15}\addtolength{\tabcolsep}{6pt}
\caption{The framework of the ADMM.\label{ADMM}}
\begin{center}
\begin{tabular}{>{\small}l}
\hline
\textbf{~~~~ADMM:~Alternating Direction Method of Multipliers}        \\\hline
~~~~Initialize $\mu,~\lambda,~\tau,~\Sigma^0,~\Lambda^0;$        \\
\textbf{~~~~Repeat}\\
~~~~~~~~~~~~~~Compute $\Gamma^{k+1}=(\Sigma^{k}+\mu\Lambda^k-\mu\tau I)_+$; \\
~~~~~~~~~~~~~~Compute $\Sigma^{k+1}=\frac{\mu}{\mu+1}\textmd{Shrink}(\Sigma_n+\frac{1}{\mu}\Gamma^{k+1}-\Lambda^k,\lambda);$ ~~~~~~ \\
~~~~~~~~~~~~~~Compute $\Lambda^{k+1}:=\Lambda^k-\frac{1}{\mu}(\Gamma^{k+1}-\Sigma^{k+1}).$ \\
\textbf{~~~~Untill} Convergence\\
\textbf{~~~~Return}\\
\hline
\end{tabular}
\end{center}
\end{table}

To end this section, we prove that the sequence $(\Gamma^k, \Sigma^k, \Lambda^k)$ produced by the alternating direction method of multipliers (Table \ref{ADMM}) converges to $(\hat{ \Gamma}, \hat{\Sigma},\hat{\Lambda})$, where $( \hat{ \Gamma},\hat{\Sigma})$ is an optimal solution of (\ref{est2}) and $\hat{\Lambda}$ is the optimal dual variable. Now we label some necessary notations for the ease of presentation. Let $H$ be a $2p\times2p$ matrix
defined as
$$H=\left(
      \begin{array}{cc}
        \frac{1}{\mu} \emph{\textbf{I}}_{p\times p} & 0 \\
        0 & \mu\emph{\textbf{I}}_{p\times p} \\
      \end{array}
    \right),
$$
the weighted norm $\left\|\cdot\right\|^{2}_{H}$ stands for $\left\|V\right\|^{2}_{H}:=\langle V,HV\rangle$ and the corresponding inner product $\langle\cdot,\cdot\rangle_{H}$ is $\langle U,V\rangle_H:=\langle U,HV\rangle$. Before presenting the main theorem with regard to the global convergence of
ADMM, we introduce the following lemma.
\begin{Lemma}\label{Conle} Assume that $(\hat{\Gamma},\hat{\Sigma})$ is an optimal solution of $\left(\ref{est2}\right)$ and $\hat{\Lambda}$ is the corresponding optimal dual variable associated with the equality constraint $\Sigma-\Gamma=0$. Then the sequence
$\left\{(\Gamma^k, \Sigma^k, \Lambda^k)\right\}$  produced by ADMM satisfies
\begin{eqnarray}\label{conle1}\|\hat{V}-V^{k}\|^{2}_{H}-\|\hat{V}-V^{k+1}\|^{2}_{H}\geq\|V^{k+1}-V^k\|^{2}_{H},\end{eqnarray}
where $V^{k}=(\Lambda^k,\Sigma^k)^{\top}$ and $\hat{V}=(\hat{\Lambda},\hat{\Sigma})^{\top}$.\end{Lemma}

\noindent Based on the lemma above, the convergent theorem can be derived immediately.
\begin{Theorem}
\label{Conth}The sequence $\left\{(\Gamma^k, \Sigma^k, \Lambda^k)\right\}$  generated by Algorithm 1 from any starting point
converges to an optimal solution of $\left(\ref{est2}\right)$.\end{Theorem}

\section{Numerical Simulations}\label{sec4}
In this section we will exploit the proposed method ADMM to tackle two examples, one of which possessed the block structured population covariance matrix, and another utilized the banded population covariance matrix. Actually as the constraint $\Sigma\succeq0$, our proposed model (\ref{est}) is equivalent to
\begin{eqnarray}\label{este}\hat{\Sigma}=\underset{\Sigma\succeq0}{\textrm{argmin}}~\frac{1}{2}\left\| \Sigma-\left(\Sigma_{n}-\tau I\right)\right\|_{F}^{2}+\lambda\left\|\Sigma\right\|_{1}.\end{eqnarray}
 So similar to the method in \cite{XMZ}, one can solve the soft-thresholding estimator $$\Sigma_{st}:=\underset{\Sigma}{\textrm{argmin}}~\frac{1}{2}\left\| \Sigma-\left(\Sigma_{n}-\tau I\right)\right\|_{F}^{2}+\lambda\left\|\Sigma\right\|_{1}$$ to initialize the $\Sigma^0$. If the derived $\Sigma_{st}\succeq0$ then the recovered sparse and low-rank semidefinite estimator $\hat{\Sigma}=\Sigma_{st}$. In our stimulation, we uniformly initialize $\Lambda^0$ as the matrix with all entries being 1, $\Sigma^0$ as  zero matrix and $\Sigma_{st}$ respectively. Unlike $\lambda$ and $\tau$, $\mu$ does not change the final covariance estimator, thus we fixed $\mu=1$ just for simplicity and the stop criteria is set as
 $$\max\left\{~\|\Gamma^{k+1}-\Gamma^{k}\|_F,~\|\Sigma^{k+1}-\Sigma^{k}\|_F~\right\}\leq 5\times10^{-4}.$$
  For the sample dimensions, we always take $n=50$ and $p=100,200,500,1000$.
\subsection{Example I: Block Structure}\label{Ex1}
Analogous to the model, modified slightly here, emerged in \cite{RSV} who synthesized  $n$ samples $X_l\sim N(0,\Sigma_0)$ for a block diagonal population covariance matrix $\Sigma_0\in \mathbb{R}^{p\times p}$, we will use $K(=5,10,20)$ blocks of random sizes, and each block is generated by $vv^{\top}$ where the entries of $v$ are drawn i.i.d. from the uniform distribution on $[-1,1]$. Evidently, the rank of $\Sigma_0$ produced in the way is $K$. Corresponding MATLAB code of generating $X_l$ is $X_l=\Sigma_0^{1/2}randn(p,1),~l=1,2,\cdots,n$, thereby deriving the sample covariance matrix
$$\Sigma_n=\frac{1}{n-1}\sum_{l=1}^{n}\left(X_{l}-\bar{X}\right)\left(X_{l}-\bar{X}\right)^{\top}.$$
What is worth mentioning is that if $\Sigma_0$ is a positive definite matrix, the solution our method obtains would be a positive definite matrix with full rank. But fortunately, compared to some largest singular values of $\Sigma_0$, the left are relatively small so that can be ignored. Here, therefore, we consider the approximate rank (\ref{apr}). In addition, we say the sparsity of a matrix $A=\left(a_{ij}\right)\in\mathbb{R}^{n\times p}$ by
$$sp(A)=\frac{\textmd{Card}\left\{~(i,j):a_{ij}\neq 0~\right\}}{n\times p}.$$
\makeatletter
          \def\@captype{figure}
          \makeatother
          \begin{center}
          \includegraphics[width=16cm]{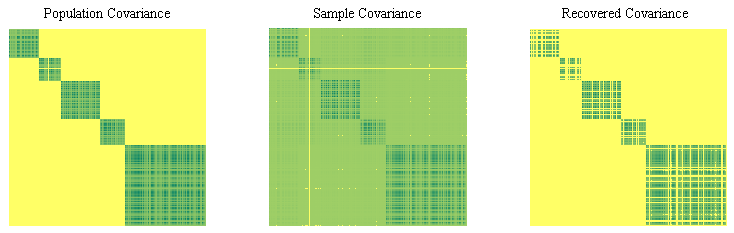}
          \end{center}
\vspace{-8mm}
          \makeatletter
          \def\@captype{figure}
          \makeatother
          \begin{center}
          \includegraphics[width=16cm]{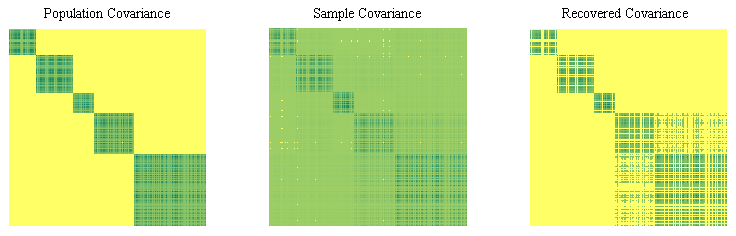}
          \end{center}
           \vspace{-7mm}
          \caption{Covariance estimation with $K=5$, $\lambda=0.15, \tau=0.2$ and $p=200$. In the row above, $FPR=0.092,TPR=1$, whilst $FPR=0.081,TPR=0.834$
          in the row below.\label{fig1}}

Apart from the approximate $ar$ and the sparsity $sp$ of the sparse and low-rank semidefinite estimator $\hat{\Sigma}=(\widehat{\sigma}_{ij})_{1\leq i,j\leq p}$, we also take advantage of other two types of errors to show the selection performance of our proposed method ADMM:
$$FPR: = \frac{\textmd{Card}\left\{(i,j):\sigma_{0ij}\neq0\&\widehat{\sigma}_{ij}=0\right\}}{\textmd{Card}~\left\{~(i,j):~\widehat{\sigma}_{ij}=0~\right\}},~~~~~~~~TPR: = \frac{\textmd{Card}\left\{(i,j):\sigma_{0ij}\neq0\&\widehat{\sigma}_{ij}\neq0\right\}}{\textmd{Card}~\left\{~(i,j):~\widehat{\sigma}_{ij}\neq0~\right\}},$$
where $FPR$ stands for the false positive rate, which means the rate of significant variables that are unselected over the whole zero entries,
 and $TPR$ denotes the true positive rate, which implies the ratio of significant variables that are selected over the entire none zero elements.
 
For more visualized purpose, we plot the Population Covariance, Sample Covariance an the Recovered Covariance. From Figures \ref{fig1} and \ref{fig2}, the left Population Covariance is $\Sigma_0$, the median Sample Covariance stands for $\Sigma_n$ and the right Recovered Covariance denotes $\hat{\Sigma}$. The yellow region stands for the sparse area in which the values are quite close (or most of them equal) to zero, while the green zone is the place where the none zero entries locate. Moreover the deeper the green color is, the larger the value stands. Evidently the recovered covariance matrices are quite dependent on the sample covariance matrix.

 \makeatletter
          \def\@captype{figure}
          \makeatother
          \begin{center}
          \includegraphics[width=16cm]{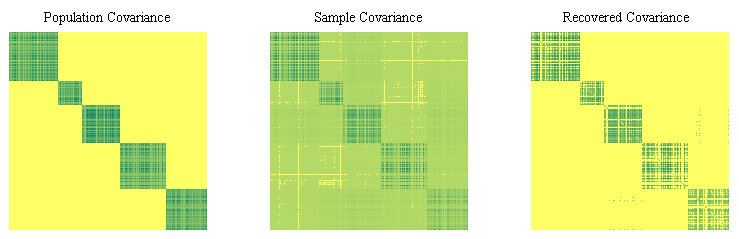}
          \end{center}
           \vspace{-8mm}
          \makeatletter
          \def\@captype{figure}
          \makeatother
          \begin{center}
          \includegraphics[width=16cm]{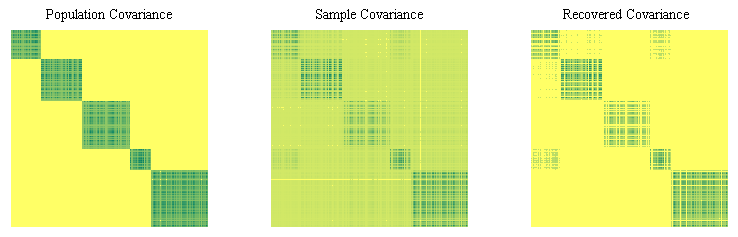}
          \end{center}
          \vspace{-7mm}
          \caption{Covariance estimation with $K=5$, $\lambda=0.15, \tau=0.2$ and $p=500$. In the row above, $FPR=0.098,TPR=0.99$, whilst $FPR=0.071,TPR=0.904$
          in the row below.\label{fig2}}

We then report average results over 100 runs. Timing (in seconds) was carried out on a CPU 2.6GHz desktop. In computation, the sparse parameter $\lambda$ and the low-rank parameter $\tau$ are given in corresponding stimulations respectively. $\Sigma^0$ ia taken as the zero matrix.  As we mentioned before the rank of the produced $\Sigma_0$ actually is given, i.e., $r(\Sigma_0)=ar(\Sigma_0)=K$. Table \ref{tab1} shows the performance of our approach under different dimensions $p=100,200,500,1000$ and various blocks $K$. Obviously, all the $ar(\hat{\Sigma})$ nearly tends to the true rank $K$. The values of $FPR$ and $TPR$ are quite desirable, which manifests that the selection performance of ADMM is very well. Moreover, the CPU time reveals the method runs relatively fast. In addition, with the $K$ rising, for example $p=500$, the $FPR$ is decreasing while $FPR$ is increasing, and the time spent by the method is also ascending. In addition, for small dimension such as $p=200$ and large block such as $K=50$, ADMM performs unstably because various stimulations can not be recovered. 
\vspace{-2mm}
 \begin{table}[h]
 \caption{The performance of ADMM under different dimensions $p$ and blocks $K$. \label{tab1}}
{\renewcommand\baselinestretch{1.1}\selectfont
\begin{center}
\begin{tabular}{c c c c c c c c }\hline
\multicolumn{8}{c}{$n=50,~~ar(\Sigma_0)=K,~~\lambda=0.25,~~\tau=0.5,~~\Sigma^{0}=0$}\\\hline
&~~$p$~~&~~$ar(\hat{\Sigma})$~~&~~$sp(\Sigma_0)$~~&~~$sp(\hat{\Sigma})$~~&~~$FPR$&$TPR$~~&~Time~~\\\hline
&100& 5.0& 0.2356  &  0.1148 &   0.1392  &  0.9863 &   0.775\\
$K=5$&200& 5.0& 0.2362 &   0.1121  &  0.1414  &  0.9911  &  3.364\\
&500&5.4&0.2225  &  0.0990  &  0.1399    &0.9748 &  32.49\\
\hline
&200& 9.90&  0.1228&    0.0565 &   0.0716  &  0.9766&    3.580\\
$K=10$&500&10.8& 0.1163  &  0.0598  &  0.0649 &   0.9371 &  33.03\\
&1000&12.0&0.1199  &  0.0604 &   0.0692 &   0.9092 & 321.6\\
\hline
&200& 19.3& 0.0601  &  0.0339 &   0.0324   & 0.8571 &   4.165\\
$K=20$&500& 20.5& 0.0569   & 0.0268   & 0.0328 &   0.9368 &  38.81\\
&1000&23.5& 0.0564  &  0.0291  &  0.0316&    0.8821 & 368.7\\
\hline
&200& 181.0&   0.0364&    0.0189 &   0.0201 &   0.8850 &  16.09\\
$K=50$&500&47.0&   0.0223 &   0.0122&    0.0127  &  0.8006  & 49.43\\
&1000&54.2&0.0221 &   0.0134 &   0.0123 &   0.7447 & 449.5\\
\hline
\end{tabular}
\end{center}
\par}
\end{table}

\vspace{-.5cm}
To simply observe the performance of our proposed method under different parameters $\lambda,\tau$ and initialized $\Sigma^0$, we fix $n=50,p=200$ and $ar(\Sigma_0)=K=5$. From Table \ref{tab11}, results of left columns of $ar(\hat{\Sigma}), FPR,TPR$ and Time are generated from the initialized $\Sigma^0=0$, and results of right columns are produced with $\Sigma^0=\Sigma_{st}$. One can easily to discern that when $\lambda=0.05$ and $\Sigma^0=\Sigma_{st}$, the performance of the method is relatively bad regardless of what $\tau$ is taken due to the $ar(\hat{\Sigma})$ and $TPR$ are undesirable. By contrast, when $\lambda=0.25$ and $0.50$, it behaves much better, particularly when  $\tau=0.25$. Moreover, with the increasing of $\lambda$, the rate $FPR$ of significant variables that are unselected is rising, even though the percentage $TPR$ of significant variables that are selected is ascending as well. By comparing the effectiveness of those two initialized point $\Sigma^0=0$ and $\Sigma_{st}$, as shown in the table, $ar(\hat{\Sigma}),FPR,TPR$ and Time generated from $\Sigma^0=0$ are basically same, which means method with zero starting point performs more stable
. But when $\lambda=0.25$ and $0.50$, ADMM with $\Sigma^0=\Sigma_{st}$ generates larger $TPR$ than that from $\Sigma^0=0$, moreover it needs less computational time in all stimulations regardless of the parameters.
 \begin{table}[h]
 \caption{The performance of ADMM under distinct parameters $\lambda,\tau$ and initialized $\Sigma^0$. \label{tab11}}
{\renewcommand\baselinestretch{1.25}\selectfont
\begin{center}
\begin{tabular}{c c c c c c }\hline
\multicolumn{6}{c}{$n=50,~~p=200,~~ar(\Sigma_0)=K=5$}\\\hline
$\lambda$&~~$\tau$~~&~~$ar(\hat{\Sigma})$~~&~~$FPR$&$TPR$~~&~Time~~\\\hline
&0.05&5.0~~12.0&0.1365~~0.0565&    1.0000~~0.5666&    3.379~~2.716\\
$0.05$&0.25&5.5~~8.5&0.1313~~0.0507&    0.8813~~0.4783&    4.857~~3.803\\
&0.50&5.0~~7.5&0.1303~~0.0615 &   0.9967~~0.4273 &   3.604~~4.040\\
\hline
&0.05&5.0~~5.0&0.1371~~0.1362 &   0.9980~~0.9843&    3.584~~1.885\\
$0.25$&0.25&5.0~5.0&0.1367~~0.1365 &   0.9972~~0.9972 &   3.963~~2.897\\
&0.50&5.2~~5.2&0.1318~~0.1317&    0.9653~~0.9661&    4.134~~3.477\\
\hline
&0.05&5.0~~5.0& 0.1375~~0.2008&    1.0000~~1.0000&    3.319~~1.822\\
$0.50$&0.25&5.4~~5.8&0.1575~~0.2169&    0.9432~~1.0000&    4.068~~2.580\\
&0.50&5.0~~6.0&0.1451~~0.2035&    0.9881~~1.0000&    3.933~~3.791\\
\hline
\end{tabular}
\end{center}
\par}
\end{table}
\vspace{-1cm}
\subsection{Example II: Banded Structure}
In this part we consider the population covariance matrix with banded structure which has been emerged in \cite{BLa,CL,XMZ}. To be more exact, the  population covariance matrix $\Sigma_0=(\sigma_{0ij})_{1\leq i,j\leq p}\in \mathbb{R}^{p\times p}$ has the following formula
$$\sigma_{0ij}=\max\Big\{~1-\frac{1}{10}|i-j|,~0~\Big\}=\Big(~1-\frac{1}{10}|i-j|~\Big)_+.$$
We first report average results over 100 replicators and take the sparse parameter $\lambda=0.5$ and the low-rank parameter $\tau=0.75$ respectively. Information listed in Table \ref{tab2} shows the performance of our approach under different dimensions $p=100,200,500,1000$ and two distinct starting point $\Sigma^{0}=0$ and $\Sigma_{st}$.
\begin{table}[h]
 \caption{The performance of ADMM over 100 simulations under different dimensions $p$ and $\Sigma^0$. \label{tab2}}
{\renewcommand\baselinestretch{1.2}\selectfont
\begin{center}
\begin{tabular}{ c c c c c c c c c}\hline
\multicolumn{8}{c}{$n=50,~~\lambda=0.5,~~\tau=0.75$}\\\hline
&~~$p$~~&~~$ar(\Sigma_0)$~~&~~$ar(\hat{\Sigma})$~~&~~$sp(\Sigma_0)$~~&~~$sp(\hat{\Sigma})$~~&~~$FPR$&$TPR$~~&~Time~~\\\hline
&100&90& 33.3&  0.1810 &   0.1209  &  0.0718&    0.9750 &   1.421\\
$\Sigma^0=0$&200&176&66.6&0.0927  &  0.0606 &   0.0351   & 0.9851 &   7.560\\
&500&487&327.0&   0.0376 &   0.0211  &  0.0174 &   0.9791 &  63.06\\
&1000&847&334.2&0.0189&0.0133&0.0067&0.9233&851.3\\
\hline
&100&90& 33.3&   0.1810  &  0.1209  &  0.0718 &   0.9750   & 1.417\\
$\Sigma^0=\Sigma_{st}$&200&176& 66.6& 0.0927 &   0.0606  &  0.0351   & 0.9851&    7.576\\
&500&487&327.0& 0.0376  &  0.0211&    0.0174 &   0.9791&   60.02\\
&1000&847&334.2&0.0189&0.0133&0.0067&0.9233&866.2\\
\hline
\end{tabular}
\end{center}
\par}
\end{table}

As we can discern in Table \ref{tab2}, compared with $\Sigma_0$, the $ar(\hat{\Sigma})$ and $sp(\hat{\Sigma})$ are relatively small, and the former ascends while the latter descends with the rise of $p$. In addition, in Example \ref{Ex1} the block structured $\Sigma_0$ whose $rank(\Sigma_0)=ar(\Sigma_0)=K$ leads to the estimator the rank of $\hat{\Sigma}$ is also close to $K$.  Being distinct with that, in this example, the $ar(\Sigma_0)$ increases with the dimension $p$ and is not low-rank, but the recovered solution $\hat{\Sigma}$ has been rendered the relatively low-rank property. The values of $FPR$ and $TPR$ are both quite desirable, which manifests that the selection performance of ADMM is very well in this example. Moreover, the CPU time reveals the method runs extremely fast as well. In addition, under such parameters Z$\lambda=0.5,\tau=0.75$, ADMM behaves nearly identically even though the starting point $\Sigma^0$ are different.

\makeatletter
          \def\@captype{figure}
          \makeatother
          \begin{center}
          \includegraphics[width=16cm]{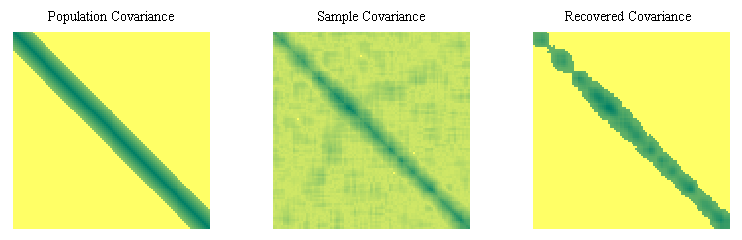}
          \end{center}
           \vspace{-8mm}
          \makeatletter
          \def\@captype{figure}
          \makeatother
          \begin{center}
          \includegraphics[width=16cm]{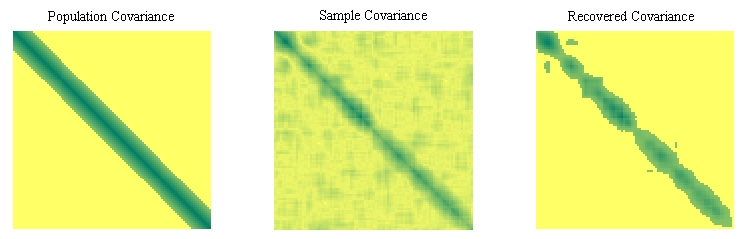}
          \end{center}
           \vspace{-2mm}
          \caption{Covariance estimation with $\lambda=0.5, \tau=0.75$ and $p=100$. In the row above, $FPR=0.062,TPR=1$, whilst $FPR=0.054,TPR=0.974$
          in the row below.\label{fig3}}
Then from Figure \ref{fig3}, one can check that the recovered covariance matrices $\hat{\Sigma}$ are quite dependent on the sample covariance matrix $\Sigma_n$, and the selection performance are relatively well because $FPR$ is pretty small while $TPR$ is close to $1$.

To simply observe the behavior of ADMM under different parameters $\lambda,\tau$ and initialized $\Sigma^0$, we fix $n=50,p=200$ and $ar(\Sigma_0)=176$. As indicated in Table \ref{tab21}, results of left columns of $ar(\hat{\Sigma}), FPR,TPR$ and Time are generated from the initialized point $\Sigma^0=0$, and results of right columns are produced with $\Sigma^0=\Sigma_{st}$. It is clear that the CPU time cost by the method with starting point $\Sigma_{st}$ is almost less than that from zero initialization. With the $\lambda$ rising, TPR generated by the method with $\Sigma^0=\Sigma_{st}$ is increasing to 1, while that with $\Sigma^0=0$ basically stabilizes at 0.97. In terms of the $ar(\hat{\Sigma})$, the proposed method with starting point $\Sigma_{st}$ will not create a lower rank solution, comparing with $\Sigma^0=0$. The reason for this phenomenon probably is that $\Sigma_{st}$ is a sparse point but not low rank; from (\ref{este}) if $\Sigma_{st}$ is an approximately semidefinite positive matrix, algorithm will stop after a few iterations which results in the solution is not a desired low-rank one.

\begin{table}[h]
 \caption{The performance of ADMM under distinct parameters $\lambda,\tau$ and initialized $\Sigma^0$. \label{tab21}}
{\renewcommand\baselinestretch{1.2}\selectfont
\begin{center}
\begin{tabular}{c c l c c  c }\hline
\multicolumn{6}{c}{$n=50,~~p=200,~~ar(\Sigma_0)=176$}\\\hline
$\lambda$&~~$\tau$~~&~~~~$ar(\hat{\Sigma})$&~~$FPR$&$TPR$~~&~Time~~\\\hline
&0.25&64.2~~126.2&  0.0356~~0.0168  &  0.9835~~0.4383 &   7.439~~5.725\\
$0.25$&0.50& 66.0~~95.2&   0.0345~~ 0.0126 &   0.9888~~0.3877&    8.682~~8.474\\
&0.75& 62.6~~78.0& 0.0381~0.0102 &   0.9780~~0.3986  &  7.852~~11.05\\
\hline
&0.25& 67.4~~103.0&0.0323~~0.0377   & 0.9528~~0.9580 &   8.996~~4.818\\
$0.50$&0.50&65.8~~80.0&0.0370~~0.0385 &   0.9834~~0.9845  &  8.445~~7.148\\
&0.75&67.8~~67.8&0.0353~~0.0353 &   0.9783~~0.9783 &8.728~~8.451\\
\hline
&0.25& 67.8~~167.0&0.0336~~0.0726 &   0.9792~~1.0000  &9.508~~2.510\\
$0.75$&0.50&64.0~~157.8& 0.0351~~0.0750 &   0.9807~~1.0000  &8.596~~1.802\\
&0.75&68.0~~166.4&0.0322~~0.0731&0.9696 ~~1.0000&8.470~~2.183\\
\hline
\end{tabular}
\end{center}
\par}
\end{table}
\section{Conclusion}
We have acquired a positive semidefinite estimator, being simultaneously sparse and low-rank, from samples of the covariance matrices through utilizing $\ell_1$ norm and nuclear norm penalties. The theoretical properties manifest that in high-dimensional settings the estimator we have constructed performs very well. Meantime, the efficient ADMM with global convergence has possessed several merits illustrated by the numerical simulations, such as less computational time and beautiful recovered effectiveness.
\section*{Acknowledgement}
The work was supported in part by the National Basic Research Program of China (2010CB732501), the National Natural Science Foundation of China (11171018, 71271021,11301022).

\begin{center}
\section*{Appendix}
\end{center}

\noindent\textbf{Proof of Lemma \ref{Conle}}~~For convenience, we denote that
$$g(\Gamma):=\tau\|\Gamma\|_*+\mathcal{I}(\Gamma\succeq0).$$
Clearly, $g:\mathbb{R}^{p\times p}\rightarrow \mathbb{R}$ is a convex function. Since  $(\hat{\Gamma},\hat{\Sigma})$ is an optimal solution of $\left(\ref{est2}\right)$, which satisfies the following KKT conditions:
\begin{eqnarray}
\label{kkt1}&&\hat{\Lambda}\in\partial\left(\tau\|\hat{\Gamma}\|_*+\mathcal{I}(\hat{\Gamma}\succeq0)\right)=\partial g(\hat{\Gamma}),\\
\label{kkt2}&&\frac{1}{\lambda}\left(-\hat{\Lambda}-\hat{\Sigma}+\Sigma_n\right)\in\partial\|\hat{\Sigma}\|_1,\\
\label{kkt3}&&\hat{\Sigma}=\hat{\Gamma}.\end{eqnarray}
~~~~Note that the optimality conditions for the first subproblem in ADMM, i.e., the
subproblem with respect to $\Gamma$ in $\left(\ref{AD1}\right)$, are given by
\begin{eqnarray}
\label{kkt11}\Lambda^k-\frac{1}{\mu}(\Gamma^{k+1}-\Sigma^{k})\in \partial g(\Gamma^{k+1}),\end{eqnarray}
this together with $\left(\ref{AD3}\right)$, i.e., $\Lambda^k=\Lambda^{k+1}+\frac{1}{\mu}\left(\Gamma^{k+1}-\Sigma^{k+1}\right)$, we have
 \begin{eqnarray}
\label{kkt12}\Lambda^{k+1}-\frac{1}{\mu}(\Sigma^{k+1}-\Sigma^{k})\in \partial g(\Gamma^{k+1}),\end{eqnarray}
Combining $\left(\ref{kkt1}\right)$ and $\left(\ref{kkt12}\right)$ and using the fact that $\partial g(\cdot)$ is a monotone operator, we get
 \begin{eqnarray}
\label{kkt13}\Big\langle\hat{\Gamma}-\Gamma^{k+1},\hat{\Lambda}-\Lambda^{k+1}+\frac{1}{\mu}(\Sigma^{k+1}-\Sigma^{k})\Big\rangle\geq0.\end{eqnarray}
~~~~The optimality conditions for the second subproblem (i.e., the subproblem with respect to $\Sigma$) in $\left(\ref{AD2}\right)$
are given by
\begin{eqnarray}
\label{kkt21}\frac{1}{\lambda}\left(-\Lambda^k-\Sigma^{k+1}+\Sigma_n-\frac{1}{\mu}(\Sigma^{k+1}-\Gamma^{k+1})\right)\in\partial\|\Sigma^{k+1}\|_1,\end{eqnarray}
this together with $\left(\ref{AD3}\right)$, i.e., $\Lambda^k=\Lambda^{k+1}+\frac{1}{\mu}\left(\Gamma^{k+1}-\Sigma^{k+1}\right)$, we have
 \begin{eqnarray}
\label{kkt22}\frac{1}{\lambda}\left(-\Lambda^{k+1}-\Sigma^{k+1}+\Sigma_n\right)\in\partial\|\Sigma^{k+1}\|_1,\end{eqnarray}
Similarly, combining $\left(\ref{kkt2}\right)$ and $\left(\ref{kkt22}\right)$, using the fact that $\partial \|\cdot\|_1$ is a monotone operator, we get
\begin{eqnarray}
\label{kkt23}\langle\hat{\Sigma}-\Sigma^{k+1},-\hat{\Lambda}+\Lambda^{k+1}-(\hat{\Sigma}-\Sigma^{k+1})\rangle\geq0,\end{eqnarray}
The summation of $\left(\ref{kkt13}\right)$ and $\left(\ref{kkt23}\right)$  gives
\begin{eqnarray}
\|\hat{\Sigma}-\Sigma^{k+1}\|_{F}^{2}&\leq& \langle\hat{\Sigma}-\Sigma^{k+1},-\hat{\Lambda}+\Lambda^{k+1}\rangle
+ \langle\hat{\Gamma}-\Gamma^{k+1},\hat{\Lambda}-\Lambda^{k+1}\rangle
+\frac{1}{\mu}\langle\hat{\Gamma}-\Gamma^{k+1},\Sigma^{k+1}-\Sigma^{k}\rangle\nonumber\\
&=&\langle\hat{\Sigma}-\Sigma^{k+1},-\hat{\Lambda}+\Lambda^{k+1}\rangle
\label{kkt4}+ \langle\hat{\Sigma}-\Sigma^{k+1}+\mu(\Lambda^{k+1}-\Lambda^{k}),\hat{\Lambda}-\Lambda^{k+1}\rangle\nonumber\\
&&+\frac{1}{\mu}\langle\hat{\Sigma}-\Sigma^{k+1}+\mu(\Lambda^{k+1}-\Lambda^{k}),\Sigma^{k+1}-\Sigma^{k}\rangle,\end{eqnarray}
where the equality because of $\Gamma^{k+1}=\Sigma^{k+1}-\mu(\Lambda^{k+1}-\Lambda^{k})$ and $\hat{\Sigma}=\hat{\Gamma}.$
Simple algebraic derivation from $\left(\ref{kkt4}\right)$ yields the following inequality:
\begin{eqnarray}
\label{kkt5}&&\|\hat{\Sigma}-\Sigma^{k+1}\|^{2}_{F}-\langle\Lambda^{k+1}-\Lambda^{k},\Sigma^{k+1}-\Sigma^{k}\rangle\nonumber\\
&\leq&\mu\langle\hat{\Lambda}-\Lambda^{k+1},\Lambda^{k+1}-\Lambda^{k}\rangle+\frac{1}{\mu} \langle\hat{\Sigma}-\Sigma^{k+1},\Sigma^{k+1}-\Sigma^{k}\rangle ,\end{eqnarray}
Rearranging the right hand side of $\left(\ref{kkt5}\right)$ using $\hat{\Lambda}-\Lambda^{k+1}=\hat{\Lambda}-\Lambda^{k}+\Lambda^{k}-\Lambda^{k+1}$ and $\hat{\Sigma}-\Sigma^{k+1}=\hat{\Sigma}-\Sigma^{k}+\Sigma^{k}-\Sigma^{k+1}$, then $\left(\ref{kkt5}\right)$ can be reduced to
\begin{eqnarray}
&&\mu\langle\hat{\Lambda}-\Lambda^{k},\Lambda^{k+1}-\Lambda^{k}\rangle+\frac{1}{\mu} \langle\hat{\Sigma}-\Sigma^{k},\Sigma^{k+1}-\Sigma^{k}\rangle\nonumber\\
&\geq&\mu\|\Lambda^{k+1}-\Lambda^{k}\|^{2}_{F}+\frac{1}{\mu}\|\Sigma^{k+1}-\Sigma^{k}\|^{2}_{F}+
\|\hat{\Sigma}-\Sigma^{k+1}\|^{2}_{F}-\langle\Lambda^{k+1}-\Lambda^{k},\Sigma^{k+1}-\Sigma^{k}\rangle ,\nonumber\end{eqnarray}
Using the notation of $V^{k}$ and $\hat{V}$, the inequality above can be rewritten as
\begin{eqnarray}
&&\langle\hat{V}-V^{k},V^{k+1}-V^{k}\rangle_{H}\nonumber\\
\label{kkt6}&\geq&\|V^{k+1}-V^{k}\|^{2}_{H}+\|\hat{\Sigma}-\Sigma^{k+1}\|^{2}_{F}-\langle\Lambda^{k+1}-\Lambda^{k},\Sigma^{k+1}-\Sigma^{k}\rangle ,\end{eqnarray}
Combining $\left(\ref{kkt6}\right)$ with the following identity
$$\|\hat{V}-V^{k+1}\|^{2}_{H}=\|V^{k+1}-V^{k}\|^{2}_{H}-2\langle\hat{V}-V^{k},V^{k+1}-V^{k}\rangle_{H}+\|\hat{V}-V^{k}\|^{2}_{H},$$
we get
\begin{eqnarray}
&&\|\hat{V}-V^{k}\|^{2}_{H}-\|\hat{V}-V^{k+1}\|^{2}_{H}\nonumber\\
&=&2\langle\hat{V}-V^{k},V^{k+1}-V^{k}\rangle_{H}-\|V^{k+1}-V^{k}\|^{2}_{H}\nonumber\\
&\geq&2\|V^{k+1}-V^{k}\|^{2}_{H}+2\|\hat{\Sigma}-\Sigma^{k+1}\|^{2}_{F}-2\langle\Sigma^{k+1}-\Sigma^{k},\Lambda^{k+1}-\Lambda^{k}\rangle -\|V^{k+1}-V^{k}\|^{2}_{H}\nonumber\\
\label{kkt7}&=&\|V^{k+1}-V^{k}\|^{2}_{H}+2\|\hat{\Sigma}-\Sigma^{k+1}\|^{2}_{F}-2\langle\Sigma^{k+1}-\Sigma^{k},\Lambda^{k+1}-\Lambda^{k}\rangle,\end{eqnarray}
Now, using $\left(\ref{kkt22}\right)$ for $k$ instead of $k+ 1$, we get,
 \begin{eqnarray}
\label{kkt222}\frac{1}{\lambda}(-\Lambda^{k}-\Sigma^{k}+\Sigma_n)\in\partial\|\Sigma^{k}\|_1,\end{eqnarray}
Combining $\left(\ref{kkt22}\right),\left(\ref{kkt222}\right)$ and using the fact that $\partial\|\cdot\|_1$ is a monotone function, we
obtain,
\begin{eqnarray}
\langle\Sigma^{k+1}-\Sigma^{k},-\Lambda^{k+1}+\Lambda^{k}-(\Sigma^{k+1}-\Sigma^{k})\rangle\geq0,\nonumber\end{eqnarray}
which means
\begin{eqnarray}
\label{kkt8}-\langle\Sigma^{k+1}-\Sigma^{k},\Lambda^{k+1}-\Lambda^{k}\rangle\geq\|\Sigma^{k+1}-\Sigma^{k}\|^{2}_{F}\geq0,\end{eqnarray}
By substituting $\left(\ref{kkt8}\right)$ into $\left(\ref{kkt7}\right)$, we get the desired result $\left(\ref{conle1}\right)$.\qed

\noindent\textbf{Proof of Theorem \ref{Conth}}~~From Lemma \ref{Conle}, we can get that\\
~~~~(a)~$\|V^{k+1}-V^k\|^{2}_{H}\rightarrow0$;\\
~~~~(b)~$\{V^k\}$ lies in a compact region;\\
~~~~(c)~$\|V^{k+1}-\hat{V}\|^{2}_{H}$ is monotonically non-increasing and thus converges.\\
Connecting with notations of $\|\cdot\|_H, V$ and $(a)$, it holds that $\Lambda^{k+1}-\Lambda^{k}\rightarrow0$ and $\Sigma^{k+1}-\Sigma^{k}\rightarrow0$, which together with $\left(\ref{AD3}\right)$ imply that $\Gamma^{k+1}-\Gamma^{k}\rightarrow0$ and $\Sigma^{k}-\Gamma^{k}\rightarrow0$. From $(b)$, $\{V^k\}$  has a subsequence $\{V^{k_j}\}$ that converges to $\check{V}=(\check{\Lambda},\check{\Sigma})^{\top}$, i.e., $\Lambda^{k_j}\rightarrow\check{\Lambda}$ and $\Sigma^{k_j}\rightarrow\check{\Sigma}$. Also we have $\Gamma^{k_j}\rightarrow\check{\Gamma}(:=\check{\Sigma})$ from $\Sigma^{k}-\Gamma^{k}\rightarrow0$. Therefore, $(\check{\Gamma},\check{\Sigma},\check{\Lambda})^{\top}$ is a limit point of $\left\{(\Gamma^k, \Sigma^k, \Lambda^k)\right\}$.

Note that $\left(\ref{kkt11}\right)$ and $\left(\ref{kkt21}\right)$ respectively imply that
\begin{eqnarray}
&&\check{\Lambda}+\frac{1}{\mu}(\check{\Gamma}-\check{\Sigma})\in \partial g(\check{\Gamma}),\nonumber\\
&&\frac{1}{\lambda}(-\check{\Lambda}-\check{\Sigma}+\Sigma_n-\frac{1}{\mu}(\check{\Sigma}-\check{\Gamma}))\in\partial\|\check{\Sigma}\|_1,\nonumber\end{eqnarray}
those together with $\check{\Gamma}=\check{\Sigma}$, it reduces that
\begin{eqnarray}
\label{kkt9}&&\check{\Lambda}\in \partial g(\check{\Gamma}),\\
\label{kkt10}&&\frac{1}{\lambda}\left(-\check{\Lambda}-\check{\Sigma}+\Sigma_n\right)\in\partial\|\check{\Sigma}\|_1,\end{eqnarray}
$\left(\ref{kkt9}\right), \left(\ref{kkt10}\right)$ and $\check{\Gamma}=\check{\Sigma}$ mean that $(\check{\Gamma},\check{\Sigma},\check{\Lambda})^{\top}$ is an optimal solution to $\left(\ref{est2}\right)$ . Therefore, we showed that any limit point of $\left\{(\Gamma^k, \Sigma^k, \Lambda^k)\right\}$ is an optimal solution to $\left(\ref{est2}\right)$.\qed

\end{document}